\documentclass[11pt]{amsart}
\usepackage{amsfonts,amsmath,amssymb}
\usepackage[all]{xy}

\begin{document}

\newcommand{\sX}{{\sf X}}
\newcommand{\sU}{{\sf U}}
\newcommand{\qX}{{\mathfrak X}}
\newcommand{\C}{{\mathbb C}}
\newcommand{\T}{{\mathbb T}}
\newcommand{\bpi}{{\boldsymbol{\Pi}}}
\newcommand{\bSigma}{{\boldsymbol{\Sigma}}}
\newcommand{\iso}{{\rm Iso}}
\newcommand{\Mod}{{\rm Mod}}
\newcommand{\Out}{{\rm Out}}
\newcommand{\op}{{\rm op}}
\newcommand{\Maps}{{\rm Maps}}
\newcommand{\Mor}{{\rm Mor}}
\newcommand{\Gps}{{\rm Gps}}
\newcommand{\Outgroups}{{\rm Outgroups}}
\newcommand{\ess}{{\mathcal S}}
\newcommand{\orb}{{\rm Orb}}
\newcommand{\Hom}{{\rm Hom}}
\newcommand{\bHom}{{\bf Hom}}
\newcommand{\cG}{{\mathcal G}}
\newcommand{\sGps}{{\sf Gps}}
\newcommand{\sHom}{{\sf Hom}}
\newcommand{\Spaces}{{\rm Spaces}}
\newcommand{\la}{{\langle}}
\newcommand{\ra}{{\rangle}}
\newcommand{\sO}{{\sf O}}
\newcommand{\cN}{{\mathcal N}}
\newcommand{\cpt}{{c}}
\newcommand{\Gpoids}{{\rm Gpoids}}
\newcommand{\sGpoids}{{\sf Gpoids}}
\newcommand{\Hot}{{\rm Hom}}
\newcommand{\SHot}{{\$ {\rm Hot}}}
\newcommand{\Func}{{\rm Func}}
\newcommand{\conj}{{\rm conj}}
\newcommand{\Sets}{{\rm Sets}}
\newcommand{\bPhi}{{\boldsymbol{\Phi}}}

\title{Categories of orbit types for proper Lie groupoids}
\author{Jack Morava}
\address{The Johns Hopkins University,
Baltimore, Maryland 21218}
\email{jack@math.jhu.edu}

\date{14 February 2014}
\begin{abstract}{It is widely understood that the quotient space of a topological
group action can have a complicated combinatorial structure, indexed somehow by the 
isotropy groups of the action [3 II \S 2.8]; but how best to record this structure
seems unclear. This sketch defines a database category of orbit types for a proper 
Lie groupoid (based on recent work [13-15] with roots in the theory of geometric 
quantization) as an attempt to capture some of this information.}\end{abstract}
\bigskip

\maketitle \bigskip

\section{Introduction and background}\bigskip

\noindent
A topological groupoid or stack [12]
\[
\xymatrix{
\sX \; := \;  s,t \; : X_1 \ar[r] \ar@<1ex>[r] & X_0 }
\]
is {\bf proper} if the map $s \times t : X_1 \to X_0 \times X_0$ is proper; 
such an object in the category of smooth manifolds and maps is a {\bf proper 
Lie groupoid}. The quotient 
\[
\sX \to \qX
\]
of $X_0$ by the equivalence relation thus defined is a Hausdorff topological 
space, sometimes called the coarse moduli space of $\sX$. \bigskip

\noindent
{\bf Examples:} \medskip

$\bullet$ Orbifolds [5]\medskip

$\bullet$ A topological {\bf transformation group}, defined by a group action
\[
G \times X \to X
\]
has an associated topological groupoid $[X/G]$ with $X_0 = X, \; X_1 = G \times X$;
I'll write $X/G$ for its quotient space. For instance
\medskip

$\bullet$ {\bf Toric varieties}, eg $G = \T^{n+1}/\T \cong \T^n$ acting on $X 
= \C P^n$ by
\[ 
(u_0,\dots, u_n) \cdot [z_0: \dots: z_n] = [u_0z_0: \dots: u_n z_n] \;,
\]
form a particularly accessible class of examples. Their quotient objects are 
polytopes: in the case above $X/G \cong \Delta^n$ is a simplex. The faces of 
the polytope define a stratification [see \S 4.1 below] of the quotient, with the 
interiors of the faces as strata. This defines an interesting poset, or category, 
associated to the groupoid: in this example it is the category of subsets of 
$\{0,\dots,n\}$ under inclusion. \bigskip

\noindent
An earlier paper [11] attempted to capture the sort of information encoded by
the face poset of a toric variety, for more general group actions. The present note
uses recent work on proper Lie groupoids [[16], cf also [1]] to propose a more general 
construction\begin{footnote}{See [8. 21] for approaches based on $\pi_1$ rather than 
$\pi_0$}\end{footnote}. \bigskip

\noindent
{\bf Acknowledgement} I am indebted to the organizers of the September 2013 Barcelona conference 
on homotopy type theory [2] for inspiration and hospitality, and for the opportunity to pursue 
these questions. I hope I will not be misunderstood by suggesting that classification problems
of the sort considered here have a deep and nontrivial history in philosophy [17].
\bigskip

\section{Some technical preliminaries} \bigskip

\noindent
{\bf 2.1 Definition} A reasonable space $X$ has a universal map $X \to \pi_0X$ 
to a discrete set, defined by the adjoint to the inclusion of the category 
of sets into that of topological spaces. The diagram
\[
\xymatrix{
\sX \ar[d] \; :  & X_1 \ar[d] \ar[r] \ar@<1ex>[r] & X_0 \ar[d] \ar@{.>}[r] & \qX \ar[d] \\
\pi_0\sX\ : & \pi_0X_1 \ar[r] \ar@<1ex>[r] & \pi_0X_0 \ar@{.>}[r] & \pi_0\qX }
\]
extends $\pi_0$ to a functor from topological to discrete groupoids, such that
\[
\pi_0[X/G] \; \cong \; [\pi_0(X)/\pi_0(G)] \;,
\]
with $\pi_0(X)/\pi_0(G) \cong \pi_0(X/G)$ for reasonable actions. The natural transformation
\[
\pi_0(\sX \times {\sf Y}) \to \pi_0 \sX \times \pi_0 {\sf Y} 
\] 
is an isomorphism in such cases. \bigskip

\noindent
{\bf 2.2} Regarding groups as categories with a single object defines a 
two-category $(\sGps)$ of groups. The set of homomorphisms from $G_0$ to 
$G_1$ has an action of $G_1$ by conjugation, defining a groupoid
\[
\Hom_\sGps(G_0,G_1) := [\Hom(G_0,G_1)/G_1^\conj] 
\]
of morphisms from $G_0$ to $G_1$.\bigskip  

\noindent
There are many variations on this theme, eg the topological two-category 
$(\sGps_\cpt)$ defined by compact groups and continuous homomorphisms. I 
will write $(\sGps^+)$ (resp. $(\sGps^+_\cpt)$) for the 
subcategories with such groups as objects, and spaces $\Hom_c^+(G_0,G_1)$ 
of continuous {\bf one-to-one} homomorphisms as maps. \bigskip

\noindent
The construction which is the identity on objects, and is the functor
\[
[\Hom_c(G_0,G_1)/G_1^\conj] \to \Hom_{\pi_0\sGps}(G_0,G_1) := \pi_0[\Hom_c(G_0,G_1)/G_1^\conj]
\] 
on morphism categories, defines a monoidal two-functor
\[
(\sGps_\cpt) \to (\pi_0\sGps_\cpt) 
\]
[and similarly for $(\sGps^+_\cpt)$]. \bigskip

\noindent 
\section{Groupoids of fixed-points with level structure} \bigskip

\noindent
{\bf 3.1 Definition:} If $\sX$ is a proper topological groupoid, and $H$ is a 
compact Lie group, let
\[
X(H)_0 := \{(x,\phi)\:|\: x \in X_0, \phi: H \to \iso(x) \in \sGps^+_\cpt \}
\]
and let $X(H)_1$ be the set of commutative diagrams of the form
\[
\xymatrix{
H \ar@{>->}[r]^{\phi'} \ar@{.>}[d]^\gamma & \iso(x') \ar[d]^{g-\conj} \\
H \ar@{>->}[r]^\phi & \iso(x) }
\]
(with $g: x' \to x \in X_1$).  The resulting proper topological groupoid 
$\sX(H)$ is a model for the subgroupoid of $\sX$ defined by points fixed by a 
group isomorphic to $H$. There is a forgetful morphism $\sX(H) \to \sX$, but 
it can't be expected to be the inclusion of a subgroupoid. \bigskip

\noindent
{\bf Proposition:} $H \mapsto \sX(H)$ defines a (two-)functor $\sX(\bullet)$ 
from $(\sGps^+_\cpt)$ to the two-category $(\sGpoids_\cpt)$ of proper 
topological groupoids. \bigskip

\noindent
{\bf Proof:} First of all, if $\xymatrix{ \alpha: H_0 \ar@{>->}[r] & H_1} \in 
(\sGps^+_\cpt)$ then 
\[ 
\xymatrix{
H_0 \ar@{>->}[r]^\alpha \ar@{.>}[d]^{\gamma^\alpha} & H_1 \ar@{>->} [r]^{\phi'_1} 
& \iso(x') \ar[d]^{g-\conj} \\
H_0 \ar@{>->}[r] & H_1 \ar@{>->}[r]^{\phi_1} & \iso(x)}
\]
defines a functor 
\[
\alpha^{H_0}_{H_1} : \sX(H_1) \to \sX(H_0) \;.
\]
Moreover, if $\alpha: H \to H$ is an inner automorphism of $H$ (ie $\alpha$ 
is conjugation by $a \in H$) then there is a natural equivalence
\[
\alpha^H_H \cong {\bf 1}_{\sX(H)}  
\]
defined by the commutative diagram
\[
\xymatrix{
H \ar[r]^\alpha \ar@{.>}[d] & H \ar@{>->}[r]^\phi & \iso(x) \ar[d]^{\phi(a^{-1})} \\
H \ar[r]^{1_H} & H \ar@{>->}[r]^\phi & \iso(x) \;.}
\]
$\Box$ \bigskip

\noindent
{\bf 3.2 Claim:} For any $\sX$ as above, there is a commutative diagram
\[
\xymatrix{
\sX \ar[r] \ar[d]^\iso & (\sGpoids_{\cpt *}) \ar[d] \\
(\sGps_\cpt) \ar[r]^{\sX[\bullet]} & (\sGpoids_\cpt) }
\]
(with the category of pointed proper groupoids in the upper right corner, and 
the forgetful map to the category of proper groupoids along the right-hand 
edge). The left-hand vertical map sends $x \in X_0$ to its isotropy group, 
and the top horizontal map sends $x$ to $\sX(\iso(x))$, with $x$ as 
distinguished point. \bigskip

\noindent
{\bf Corollary} The universal property of a fiber product defines a continuous 
functor
\[
\sX \to \Phi_0(\sX)
\]
to the category defined by the pullback 
\[
\xymatrix{
\Phi_0(\sX) \ar@{.>}[d] \ar@{.>}[r] & (\sGpoids_{\cpt *}) \ar[d] 
\ar[r]^{\pi_0} & (\sGpoids)_* \ar[d]\\
(\sGps_\cpt) \ar[r]^{\sX[\bullet]} & (\sGpoids_\cpt) \ar[r]^{\pi_0} & (\sGpoids) }
\]  
(where the two right vertical arrows are the obvious forgetful 
functors). $\Box$ \bigskip

\noindent
It's natural to think of $\Phi_0(\sX)$ as a {\bf database} category 
[11, 18]. However I don't know how to characterize $\Phi_0$ by some universal 
property (such as being an adjoint). \bigskip

\noindent
{\bf 3.3 Example:} A proper transformation group $[X/G]$ defines a functor
\[
G > H \mapsto X^H = \{ x \in X \:|\: \iso(x) \subset H \}
\]
from the topological category $(G-\orb)$ [with closed subgroups of $G$ as 
objects, and 
\[
\Mor_{G-\orb}(H_0,H_1) = \Maps_G(G/H_0,G/H_1) =  \{g \in G\:|\: gH_0g^{-1} \subset 
H_1 \}/H_1^\conj 
\]
as morphism objects [3 I \S 10], to spaces. \bigskip

\noindent
This extends to a functor
\[
S^0[X^\bullet] : (G-\Spaces) \ni X \mapsto S^0[X^H] \in \Func(G-\orb,S^0-\Mod)
\]
which provides a model [6 V \S 9, 10, 20] for the $G$-equivariant stable category in terms of 
sheaves of spectra (ie $S^0$-modules) over $(G-\orb)$.  \bigskip

\noindent
The commutative diagram
\[
\xymatrix{
[X/G] \ar@{.>}[r] \ar@{.>}[d]^\iso & (\Sets)_* \ar[d] \\
(G-\orb) \ar[r]^{\pi_0X[*]} & (\Sets) }
\]
defines a functor from $[X/G]$ to a fiber product category $\bPhi_0[X/G]$ (with objects,
pairs consisting of a subgroup $H$ of $G$, and a component of $X^H$), analogous to the 
construction in the previous paragraph [11 \S 2.2]. The sheaf $S^\infty X^\bullet$ of spectra pulls
back to a sheaf of spectra over $\bPhi_0(X)$. 
\bigskip

\noindent
One might hope for an unstable version of this construction, applicable in the theory of
$\infty$-categories (cf eg [7 \S 5.5.6.18]); but because it depends on a presentation of
$[X/G]$ as a global quotient, it does not seem to be homotopy-invariant. 
\bigskip

\section{Proper {\bf Lie} groupoids, after Pflaum {\it et al}}\bigskip

\noindent
{\bf 4.1} A {\bf stratification} $\ess$ of a (paracompact, second countable) topological space $X$
assigns to each $x \in X$, the germ of a closed subset $\ess_x$ (containing $x$) of $X$. With suitably
defined morphisms [9 \S 1.8, 16 \S 1], stratified spaces form a category. A stratification defines 
a locally finite partition
\[
X = \coprod_{S \in \Sigma(\ess)} X_S
\]
of $X$ into locally closed subsets (called its strata), such that if $x \in X_S$ then $\ess_x$ is the
associated set germ. 
\bigskip

\noindent
Very interesting recent work of M. Pflaum {\it et al}  [building 
on earlier work of Weinstein and Zung ([22]; cf also [13]) shows that 
\bigskip

\noindent
{\bf Theorem} [16 Theorem 5.3, Cor 5.4] The quotient space $\qX$ of a proper Lie 
groupoid $\sX$ has a canonical {\bf Whitney} stratification. The associated 
decomposition of $X_0$ into locally closed submanifolds 
\[
X_{0(H)} \; = \; \{x \in X_0 \:|\: \iso(x) \cong H \}
\]
is indexed [16 Theorem 5.7] by (isomorphism classes of) compact Lie groups $H$. \bigskip

\noindent
{\bf 4.2 Definition}  The {\bf orbit} groupoid $\sO(x) \subset \sX$ of $x \in X_0$
\[
\sO_0(x) = \{ y \in X_0 \:|\: \exists g: y \to x \; \in X_1 \} n
\]
\[
\sO_1(x) = \{ g \in X_1 \:|\: s(g),t(g) \in \sO_0(x) \}
\]
reduces, in the case of a transformation groupoid $[X/G]$, to the groupoid
\[
[(G/\iso(x))/G] \equiv [*/\iso(x)] \;.
\]

\noindent
A {\bf slice} at $x \in X_0$ is (very roughly [13 \S 3.3-4,3.8-9]) the germ of an
$\iso(x)$-invariant submanifold of $X_0$ containing $x$, transverse to $\sO_0(x)$; for
a transformation group it is something like the image of an exponential map
\[
[\cN_x/\iso(x)] \equiv [(\cN_x \times_{\iso(x)} G)/G] \to [X/G]
\]
(where $\cN_x \in (\iso(x) - \Mod)$ is the linear representation
\[
0 \to T_x G \to T_x X_0 \to \cN_x \to 0 \;.
\]
defining the normal bundle to the orbit of $x$).\bigskip

\noindent
{\bf Theorem} [16 \S 3.11] There is an (essentially unique) slice at every object
of a proper Lie groupoid $\sX$; the corresponding set germs define the canonical stratification
[16 \S 5.4] of $\sX$.\bigskip

\noindent
{\bf Definition} The {\bf normal orbit type} of $x \in X_0$ is the equivalence class of its normal
$\iso(x)$-representation $\cN_x$. More precisely, $x_0 \sim x_1$ if there are isomorphisms
\[
\phi: \iso(x_0) \to \iso(x_1), \; \Phi: \cN_{x_0} \to \phi^*(\cN_{x_1})
\]
of groups and representations. The connected components $\nu \in \pi_0(X_0)$ of the normal orbit types 
of $\sX$ are [16 \S 5.7] the strata of the canonical partition of $X_0$.\bigskip

\noindent
The {\bf condition of the frontier} [16 Prop 5.15] asserts that if $\nu' \cap \overline{\nu} \neq 
\emptyset$ then $\overline{\nu} \supset \nu'$. This implies the existence of a partial order ($\nu > 
\nu'$) on the set $\Sigma(\sX)$ of connected components of normal orbit types for $\sX$, which can 
thus be regarded as the objects of a category [4 II \S 2.8]. Thus
\[                                                                                                                         \overline{X_{0(K)}} = \coprod_{\pi_0(X_{0(K)}) \ni \nu > \nu'} \nu' \;,
\]
\medskip

\noindent
This gives us some control of the functor $\Phi_0$ on proper Lie groupoids: \bigskip

\noindent
{\bf 4.3 Proposition} For a proper Lie groupoid $\sX$, we have isomorphisms
\[
\xymatrix{
\bigcup_{H < K} \; X_{0(K)} \times \Hom^+_c(H,K) \ar[r]^<<<<\cong & X(H)_0} 
\]
\[
\xymatrix{
\bigcup_{H < K, x \in X_{0(K)}} \sO_1(x) \times K \ar[r]^<<<<\cong & X(H)_1}
\]
and consequently 
\[
\xymatrix{
\bigcup_{H < K} \; \qX_K \times \Hom^+_c(H,K)/K^\conj \ar[r]^<<<<\cong & \qX(H)} 
\] 
(where $\qX_K \subset \qX$ is the space of orbits with isotropy group isomorphic 
to $K$). \bigskip

\noindent
{\bf 4.4 Closing remarks} \bigskip

\noindent
i) When $\sX = [X/G]$ this all simplifies a little. In particular, since
\[
X^H = \bigcup_{H < K <G} X_K \;,
\]
$\bPhi_0[X/G]$ is essentially just the quotient of $\Phi_0[X/G]$ which 
collapses the morphism spaces $\Hom_c^+(K,H)/K^\conj$.\bigskip

\noindent
ii) The subspaces $X_{0(K)}$ are disjoint unions of strata $\nu$ indexed by slice representations
\[
K \to {\rm Aut}(\cN_\nu) \;.
\]
The resulting family of vector spaces over $\sX(H)$ pulls back to a fibered category 
\[
\cN(\sX) \to \Phi_0(\sX) \;.
\]
This seems to provide a natural repository for Noether's theorem (which associates
conserved quantities to elements of the Lie algebra of symmetries of states of a physical
system) [11 \S 4.1]. \bigskip

\noindent
iii) I don't know how generally one can associate a stratification to a topological groupoid.
There are many interesting examples, coming from locally compact groupoids (eg the Thom-Boardman theory
of singularities of smooth maps [18]), or from infinite-dimensional examples (Ebin's category of 
Riemannian metrics up to diffeomorphism, Vassiliev's finite-type invariants of immersions, \dots), 
where a more general theory would be very interesting. The existence and good behavior of slices seem 
to be an essential requirement for such a theory. 

\newpage

\bibliographystyle{amsplain}

\end{document}